\documentclass[10pt]{amsart}

\usepackage{algorithm}
\usepackage{algpseudocode}
\usepackage{amsmath}
\usepackage{amssymb}
\usepackage{amsthm}
\usepackage{amsfonts}
\usepackage{caption}
\usepackage{hyperref}
\usepackage{mathtools}
\usepackage{pgfplots}
\usepackage{verbatim}

\usepackage{tikz}
\usetikzlibrary{arrows,backgrounds,calc}

\usepackage[backend=bibtex,style=alphabetic,citestyle=alphabetic]{biblatex}

\newtheorem{theorem}{Theorem}[section]

\newtheorem{conjecture}[theorem]{Conjecture}
\newtheorem{counterexample}[theorem]{Counter-example}

\newtheorem{proposition}[theorem]{Proposition}

\numberwithin{equation}{theorem}

\theoremstyle{definition}
\newtheorem{definition}[theorem]{Definition}
\newcommand{\definedterm}[1]{\textbf{#1}}

\providecommand\given{}
\newcommand\SetSymbol[1][]{%
  \nonscript\:#1\,\vert\,\allowbreak\nonscript\:\mathopen{}}
\DeclarePairedDelimiterX\set[1]\{\}{%
  \renewcommand\given{\SetSymbol[\delimsize]}
  #1
}
\DeclarePairedDelimiterX\angset[1]\langle\rangle{%
  \renewcommand\given{\SetSymbol[\delimsize]}
  #1
}

\DeclareMathOperator{\cheegerchernsimons}{\widehat{c}}
\DeclareMathOperator{\chernsimons}{CS}
\DeclareMathOperator{\complexvol}{Vol_{\mathbb{C}}}
\DeclareMathOperator{\im}{Im}
\DeclareMathOperator{\tr}{tr} \DeclareMathOperator{\vol}{vol}

\newcommand{\Isom}{\operatorname{Isom}}
\newcommand{\PSL}{\operatorname{PSL}}
\newcommand{\SL}{\operatorname{SL}}
\newcommand{\bloch}[1]{\mathcal{B}(#1)}
\newcommand{\codebit}[1]{\texttt{#1}}
\newcommand{\crossratio}[4]{[#1\mathop{:}#2\mathop{:}#3\mathop{:}#4]}
\newcommand{\dx}[1]{\, \mathrm{d} {#1}}
\newcommand{\invtracefield}[1]{\mathbb{Q}(\tr #1^{(2)})}
\newcommand{\prebloch}[1]{\mathcal{P}(#1)}
\newcommand{\prog}[1]{\texttt{#1}}
\newcommand{\tracefield}[1]{\mathbb{Q}(\tr #1)}

\begin{document}

\title[Investigation of Neumann's Conjecture]{An experimental
  investigation of Neumann's conjecture}

\author[Gilles]{Stephen Gilles} \address{Mathematics Department \\
  University of Maryland \\ 4176 Campus Dr. \\ College Park, MD 20742}
\email{sgilles@math.umd.edu}

\author[Huston]{Peter Huston}
\email{peter.huston@otterbein.edu}

\thanks{This work was supported in part by the National Science
  Foundation through the Research Experiences for Undergraduates
  (DMS\# 1359307) at the University of Maryland, College Park.}

\date{2016-09-17}


\begin{abstract}
  We use a large census of hyperbolic 3-manifolds to experimentally
  investigate a conjecture of Neumann regarding the Bloch Group. We
  present an augmented census including, for feasible invariant trace
  fields, explicit manifolds (associated to that field) that appear to
  generate the Bloch group of that field. We also make use of Ptolemy
  coordinates to compute ``exotic volumes'' of representations, and
  attempt to realize these volumes as linear combinations of generator
  volumes.  We thus present a large body of empirical support for
  Neumann's conjecture.
\end{abstract}

\maketitle

\tableofcontents

\section{Introduction}

Our primary result, described in section~\ref{sec:main-result}, is
a large census of hyperbolic 3-manifolds, organized by invariant
trace field, together with tables expressing volumes as linear
combinations of other manifolds' volumes. This census is available at
\url{http://www.curve.unhyperbolic.org/linComb}.  Examining this census
has allowed us, in many cases, to explicitly list manifolds which give
rise to elements in the Bloch groups of certain invariant trace fields
that span all the elements of that Bloch group we observed, supporting
a conjecture of Neumann.

\subsection{Neumann's conjecture}

For a hyperbolic $3$-manifold $M$, we denote the element it induces in
the Bloch Group (see section~\ref{subsec:bloch-group}) as $[M]$. The
conjecture we primarily studied was given in \parencite{neumann06},
and stated in this form by \parencite{garoufalidis2015}:

\begin{conjecture}[Neumann]
  \label{conj:neumann-main}
  Let $F \subset \mathbb{C}$ be a number field not contained in
  $\mathbb{R}$, and let $\mathcal{M}$ be the set of manifolds with
  invariant trace fields contained in $F$, \[ \mathcal{M} := \set{
  M \cong \mathbb{H}^3/\Gamma \given
      \invtracefield{\Gamma} \subseteq F }. \] Let
  $\mathcal{N} = \set{[M] \given M \in \mathcal{M}}$ be elements of the
  Bloch group determined by manifolds in $\mathcal{M}$ (see
  section~\ref{subsec:bloch-group}). Then integral
  combinations of elements of $\mathcal{N}$ generate $\bloch{F}$.
\end{conjecture}

\subsection{Neumann's weaker conjecture}

Neumann has also proposed a weaker conjecture which our results also
support.

\begin{conjecture}[{\parencite[Conjecture~1]{neumann06}}]
  \label{conj:neumann-weaker}
  Every non-real concrete number field $k$ arises as the invariant trace
  field of some hyperbolic manifold.
\end{conjecture}

In full form, this conjecture also addresses quaternion algebras, but
we only examine invariant trace fields. Relevant results are discussed
in section~\ref{sec:manifold-existence}.

\subsection{A stronger, false conjecture}

We also considered whether conjecture~\ref{conj:neumann-main} could be
strengthened somewhat to the following form, which directly regards
volumes of manifolds rather than Bloch group elements.

\begin{conjecture}
  \label{conj:false-strengthening}
  Let $F \subset \mathbb{C}$ be a number field not contained in
  $\mathbb{R}$.

  Let $S = \set{ \vol(M) \given M \text{ a hyperbolic $3$-manifold
      with invariant trace field } F}$. The lattice generated by $S$
  is linearly spanned by some $\set{\vol(N_1), \dotsc,
  \vol(N_{r_2})} \subset S$, where $r_2$ is the number of complex places
  of $F$ (and each $N_i$ is a hyperbolic 3-manifold).
\end{conjecture}

In section~\ref{sec:strengthening-counterexample}, we provide a
counterexample to conjecture~\ref{conj:false-strengthening}.

\subsection{Exotic volumes}

Recall that for a closed hyperbolic 3-manifold and a representation
$\rho : \pi_1(M) \to \SL(n, \mathbb{C})$, Cheeger-Chern-Simons
invariant $\cheegerchernsimons(\rho)$ is given by
\begin{equation}
  \label{eq:cheeger-chern-simons}
  \cheegerchernsimons(\rho) = \frac{1}{2} \int_{M} s^{\ast} \left( \tr \left( A \wedge \dx{A} + \frac{2}{3} A \wedge A \wedge A \right) \right) \in \mathbb{C} / 4\pi^2\mathbb{Z},
\end{equation}
with $E_{\rho}$ the flat $\SL(n, \mathbb{C})$-bundle with holonomy
$\rho$, with $A$ the flat connection in $E_{\rho}$, and with $s$ a
section of $E_{\rho}$.

\begin{definition}
  For a representation $\rho: \pi_1(M) \to \SL(n, \mathbb{C})$, the
  \definedterm{complex volume} $\complexvol(\rho)$ of $\rho$ is
  \[ \complexvol(\rho) = i \cheegerchernsimons(\rho). \] If $\rho$ is
  $\rho_{\text{geo}}$, the geometric representation of a hyperbolic
  $3$-manifold $M$, then
  \[ \complexvol(\rho_{\text{geo}}) = \vol(M) + i \chernsimons(M), \]
  for $\chernsimons(M)$ the Chern-Simons invariant of $M$.

  When $\rho$ is not $\rho_{\text{geo}}$, we call the real part of
  $\complexvol(\rho)$ an \definedterm{exotic volume} of $M$, to
  distinguish it from the geometric volume.
\end{definition}

By means of the Ptolemy coordinates of \parencite{garoufalidis2015},
we were able to compute certain exotic volumes of manifolds to at
least 50 decimal places of precision, and to express them as linear
combinations of geometric volumes of manifolds associated to the same
invariant trace field.

\subsection*{Acknowledgments}

The authors would like to thank Christian Zickert for guidance and
advice, and the 2014 MAPS-REU program at the University of Maryland,
led by Kasso Okoudjou, as part of which this work was undertaken.

\section{Background}

\subsection{Cross-ratios}

It is a feature of hyperbolic geometry that the volume of a complete
3-manifold (if finite) is actually a topological invariant (see e.g.\
\parencite[83]{benedetti1992lectures} for a proof). For an ideal
hyperbolic 3-simplex $\Delta = (a,b,c,d)$ defined by four points in
$\partial \mathbb{H}^3 \cong \mathbb{C} \cup \set{\infty}$,
$\vol(\Delta)$ can be determined knowing just the cross-ratio
$\crossratio{a}{b}{c}{d} = \frac{(c-a)(d-b)}{(d-a)(c-b)}$.  Since
$\Isom^{+}(\mathbb{H}^3) \cong \PSL(2, \mathbb{C})$ is 3-transitive,
for the purposes of volume calculation we may assume $\Delta$ to be of
the form $(\infty,0,1,z)$, which conveniently has cross-ratio $z$. The
volume of $\Delta$ is given by
\begin{equation}
  \label{eq:D}
  D(\Delta) = \im \left( \int_0^1 \frac{\log (1 - tz)}{t} \dx{t} \right) + \arg(1 - z) \log |z|.
\end{equation}

Given a triangulation of a hyperbolic 3-manifold, one may compute
cross-ratios using gluing equations \ref{eq:cusp} and
\ref{eq:edge}. Using cross-ratios, one may compute volumes for each
simplex in the triangulation; we refer to
section~\ref{subsec:gluing-equations} for more details. Summing those
volumes gives the volume of the manifold: for $\set{\Delta_i}$ the
simplices in a triangulation of $M$,
\begin{equation}
  \label{eq:sum-of-simplex-vols}
  \vol(M) = \sum_i D(\Delta_i).
\end{equation}

In our work, we used experimental data to examine relations between
the volumes of hyperbolic 3-manifolds which share the same invariant
trace field.  We begin by briefly reviewing the basics of the algebraic
invariants of hyperbolic 3-manifolds relevant to our study.

\subsection{The invariant trace field}

\begin{definition}
  For a hyperbolic manifold $M = \mathbb{H}^3 / \Gamma$ defined by a
  discrete subgroup $\Gamma$ of $\PSL(2, \mathbb{C})$, the
  \definedterm{trace field} of $\Gamma$, denoted
  $\tracefield{\Gamma}$, is
  \[ \mathbb{Q}( \set{ \tr \gamma : \pi(\gamma) \in \Gamma } ) \]
  where $\pi: \SL(2, \mathbb{C}) \to \PSL(2, \mathbb{C})$ is the
  standard projection map.
\end{definition}

Since for any matrices $a, b \in \SL(2, \mathbb{C})$ we have
$\tr (aba^{-1}) = \tr(b)$, the trace field of $\Gamma$ is a
conjugation invariant of $\Gamma$. Unfortunately, it is not a
commensurability invariant, since it is possible to create finite
degree extensions of $\Gamma$ which extend the trace field: we refer
to \parencite[116]{maclachlanreid} for an explicit example. A slight
modification smooths over this difficulty.

\begin{definition}
  For a non-elementary subgroup $\Gamma \subset \PSL(2, \mathbb{C})$,
  the \definedterm{invariant trace field} of $\Gamma$ is the trace
  field $\invtracefield{\Gamma}$, where
  \[ \Gamma^{(2)} := \angset{ \gamma^2 \given \gamma \in \Gamma }. \]
\end{definition}

\subsection{The Bloch group}
\label{subsec:bloch-group}


\begin{definition}
  For any field $F$, the \definedterm{Pre-Bloch group on $F$}, denoted
  $\prebloch{F}$, is defined as
    \[ \prebloch{F} = \frac{\langle \set{z \in F
          \given z \ne 0_F, 1_F} \rangle}{[z_1] - [z_2] +
        \left[\frac{z_2}{z_1}\right] - \left[\frac{1 - z_2}{1 -
            z_1}\right] + \left[ \frac{1 - z_2^{-1}}{1 -
            z_1^{-1}}\right]}. \]
  \end{definition}

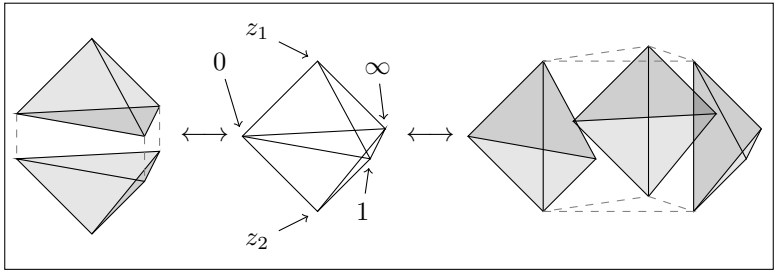
\begin{figure}[h]
  \tikzset{ c/.style={every coordinate/.try} }
  \begin{tikzpicture}[framed]
    \coordinate (A) at (0,0);
    \coordinate (Albl) at (-0.3,1);
    \coordinate (B) at (1.7,-0.3);
    \coordinate (Blbl) at (1.6,-1);
    \coordinate (C) at (1.9,0.1);
    \coordinate (Clbl) at (1.8,0.9);
    \coordinate (D) at (1,1);
    \coordinate (Dlbl) at (0.2,1.4);
    \coordinate (E) at (1,-1);
    \coordinate (Elbl) at (0.2,-1.4);

    \begin{scope}[every coordinate/.style={shift={(0,0.3)}}]
      \filldraw[opacity=0.1] ([c]A) -- ([c]B) -- ([c]C) -- cycle;
      \filldraw[opacity=0.1] ([c]A) -- ([c]B) -- ([c]D) -- cycle;
      \filldraw[opacity=0.1] ([c]D) -- ([c]B) -- ([c]C) -- cycle;
      \draw ([c]A) -- ([c]B) -- ([c]C) -- cycle;
      \draw ([c]A) -- ([c]D) -- ([c]B);
      \draw ([c]D) -- ([c]C);
    \end{scope}

    \begin{scope}[every coordinate/.style={shift={(0,-0.3)}}]
      \filldraw[opacity=0.1] ([c]A) -- ([c]B) -- ([c]C) -- cycle;
      \filldraw[opacity=0.1] ([c]A) -- ([c]B) -- ([c]E) -- cycle;
      \filldraw[opacity=0.1] ([c]E) -- ([c]B) -- ([c]C) -- cycle;
      \draw ([c]A) -- ([c]B) -- ([c]C) -- cycle;
      \draw ([c]A) -- ([c]E) -- ([c]B);
      \draw ([c]E) -- ([c]C);
    \end{scope}

    \draw[gray, dashed] ([shift={(0,0.3)}]A) -- ([shift={(0,-0.3)}]A);
    \draw[gray, dashed] ([shift={(0,0.3)}]B) -- ([shift={(0,-0.3)}]B);
    \draw[gray, dashed] ([shift={(0,0.3)}]C) -- ([shift={(0,-0.3)}]C);

    \node at (2.5,0) {$\longleftrightarrow$};

    \begin{scope}[every coordinate/.style={shift={(3,0)}}]
      \draw ([c]A) -- ([c]D) -- ([c]B) -- ([c]C) -- ([c]E) -- cycle;
      \draw ([c]A) -- ([c]B) -- ([c]E);
      \draw ([c]A) -- ([c]C) -- ([c]D);

      \node (abubble) at ([c]A) {};
      \node (alabel) at ([c]Albl) {$0$};
      \draw [->] (alabel) -- (abubble);

      \node (bbubble) at ([c]B) {};
      \node (blabel) at ([c]Blbl) {$1$};
      \draw [->] (blabel) -- (bbubble);

      \node (cbubble) at ([c]C) {};
      \node (clabel) at ([c]Clbl) {$\infty$};
      \draw [->] (clabel) -- (cbubble);

      \node (dbubble) at ([c]D) {};
      \node (dlabel) at ([c]Dlbl) {$z_1$};
      \draw [->] (dlabel) -- (dbubble);

      \node (ebubble) at ([c]E) {};
      \node (elabel) at ([c]Elbl) {$z_2$};
      \draw [->] (elabel) -- (ebubble);
    \end{scope}

    \node at (5.5,0) {$\longleftrightarrow$};

    \begin{scope}[every coordinate/.style={shift={(6,0)}}]
      \filldraw[opacity=0.1] ([c]D) -- ([c]A) -- ([c]B) -- cycle;
      \filldraw[opacity=0.1] ([c]D) -- ([c]E) -- ([c]B) -- cycle;
      \filldraw[opacity=0.1] ([c]D) -- ([c]A) -- ([c]E) -- cycle;
      \draw ([c]A) -- ([c]D) -- ([c]B) -- cycle;
      \draw ([c]A) -- ([c]E) -- ([c]B);
      \draw ([c]E) -- ([c]D);
    \end{scope}

    \begin{scope}[every coordinate/.style={shift={(8,0)}}]
      \filldraw[opacity=0.1] ([c]D) -- ([c]C) -- ([c]B) -- cycle;
      \filldraw[opacity=0.1] ([c]D) -- ([c]B) -- ([c]E) -- cycle;
      \filldraw[opacity=0.1] ([c]E) -- ([c]D) -- ([c]C) -- cycle;
      \draw ([c]D) -- ([c]C) -- ([c]B) -- cycle;
      \draw ([c]D) -- ([c]E) -- ([c]B);
      \draw ([c]E) -- ([c]C);
    \end{scope}

    \begin{scope}[every coordinate/.style={shift={(7.4,0.2)}}]
      \filldraw[opacity=0.1] ([c]A) -- ([c]D) -- ([c]C) -- cycle;
      \filldraw[opacity=0.1] ([c]A) -- ([c]D) -- ([c]C) -- cycle;
      \filldraw[opacity=0.1] ([c]A) -- ([c]E) -- ([c]C) -- cycle;
      \draw ([c]A) -- ([c]D) -- ([c]C) -- cycle;
      \draw ([c]A) -- ([c]E) -- ([c]C);
      \draw ([c]E) -- ([c]D);
    \end{scope}

    \draw[gray, dashed] ([shift={(6,0)}]D) -- ([shift={(8,0)}]D)
                        -- ([shift={(7.4,0.2)}]D) -- cycle;
    \draw[gray, dashed] ([shift={(6,0)}]E) -- ([shift={(8,0)}]E)
                        -- ([shift={(7.4,0.2)}]E) -- cycle;
  \end{tikzpicture}
  \caption{The Pachner $2$-$3$ move.}
  \label{fig:pachner-2-3}
\end{figure}

  The relation in the denominator is known as the \definedterm{Five
    Term Relation}. The variables $z_1$ and $z_2$ are to be understood
  as cross-ratios of hyperbolic simplices. The relation encodes that
  the union of two simplices with volumes $[z_1]$ and $[z_2]$ is also
  the union of three simplices with volumes $[\frac{z_2}{z_1}]$,
  $[\frac{1 - z_2}{1 - z_1}]$, and
  $[ \frac{1 - z_2^{-1}}{1 - z_1^{-1}}]$. This equivalence is called
  the Pachner $2$-$3$ move (see figure~\ref{fig:pachner-2-3}).

\begin{definition}
  For any field $F$, the \definedterm{Bloch group on $F$}, denoted
  $\bloch{F}$, is the kernel of the map
  $d : \prebloch{F} \to F^{\ast} \wedge_{\mathbb{Z}} F^{\ast}$ defined
  by \[ d: [z] \mapsto z \wedge (1-z). \]
\end{definition}

The map $d$ is an analogue of the Dehn invariant map
\parencite{neumannhilberts3rd}, which (for $F = \mathbb{C}$) fits into
an exact sequence of scissors congruence described in
e.g.\ \parencite{dupontscissors}. Thus, understanding
$\bloch{F} = \ker d$ provides insight into scissors congruence groups
and generalizations of Hilbert's third problem.

\subsection{Gluing Equations}
\label{subsec:gluing-equations}

Given a hyperbolic 3-manifold $M$ with cusps, the standard procedure
for determining $\vol(M)$ is to first triangulate $M$ into a
collection of simplices $\set{\Delta_i}$ (together with face
pairings), next to compute the cross ratio $z_i$ of each $\Delta_i$,
and finally to compute $\vol(M)$ by
equation~\ref{eq:sum-of-simplex-vols}. The first step is
well-understood: the manifolds included in the censuses of
\prog{SnapPy} include triangulations. The third step is a matter of
numerical approximation, but the second step deserves more
explanation. We present an overview of the \definedterm{gluing
  equations} here, and refer to \parencite{neumann-zagier} for a more
detailed description, including a treatment of gluing equations for
closed manifolds obtained by Dehn surgery.

Given a triangulation $\set{\Delta_i}$ of $M$, the standard method of
finding $\set{z_i}$ is to exploit the torus boundary components and
local Euclidean nature of the manifold to produce and solve a system
of equations in $\set{z_i}$. Equations~\ref{eq:cusp}, which arise from
torus boundary components, are the \definedterm{cusp equations}, and
equations~\ref{eq:edge}, which arise from local Euclideanness, are the
\definedterm{edge equations}. Together, they are the gluing equations.

\begin{figure}[h]
  \begin{tikzpicture}[framed]
    \coordinate (A) at (0.4,0);
    \coordinate (B) at (1.4,-0.6);
    \coordinate (C) at (2,0);
    \coordinate (D) at (2,2);
    \coordinate (A') at ($(A)!0.6!(D)$);
    \coordinate (B') at ($(B)!0.6!(D)$);
    \coordinate (C') at ($(C)!0.6!(D)$);

    \coordinate (E) at (4,-0.6);
    \coordinate (F) at (2.5,-0.6);
    \coordinate (G) at (3,0);
    \coordinate (H) at (3,2);
    \coordinate (E') at ($(E)!0.6!(H)$);
    \coordinate (F') at ($(F)!0.6!(H)$);
    \coordinate (G') at ($(G)!0.6!(H)$);

    \draw[dashed] (C) -- (A);
    \draw (A) -- (B) -- (C);
    \filldraw[opacity=0.1] (A') -- (B') -- (C') -- cycle;
    \draw (A') -- (B') -- (C') -- cycle;
    \draw (A) -- (A');
    \draw (B) -- (B');
    \draw (C) -- (C');
    \draw[dashed] (A') -- (D);
    \draw[dashed] (B') -- (D);
    \draw[dashed] (C') -- (D);

    \draw[dashed] (F) -- (G) -- (E);
    \draw (E) -- (F);
    \filldraw[opacity=0.1] (E') -- (F') -- (G') -- cycle;
    \draw (E') -- (F') -- (G') -- cycle;
    \draw (E) -- (E');
    \draw (F) -- (F');
    \draw[dashed] (G) -- (G');
    \draw[dashed] (E') -- (H);
    \draw[dashed] (F') -- (H);
    \draw[dashed] (G') -- (H);

    \draw[dotted] (B) -- (F);
    \draw[dotted] (B') -- (F');
    \draw[dotted] (C) -- (G);
    \draw[dotted] (C') -- (G');
    \draw[dotted] (D) -- (H);
  \end{tikzpicture}
  \caption{Piecing together a torus from a triangulation.}
  \label{fig:constructing_torus}
\end{figure}
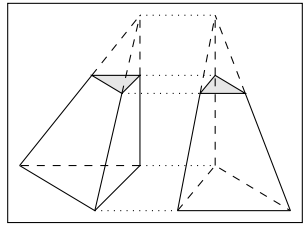

Since, in the triangulation of $M$, a link component appears at the
vertices of each simplex, a torus can be determined by truncating the
simplices, then making edge identifications to match face
identifications, as in figure~\ref{fig:constructing_torus}.

Since the meridian and the longitude of a torus are homotopically circles,
it should be possible to require that any path following a meridian or
longitude winds exactly once. In other words, it should be possible to
construct a system of equations, with variables representing various
angles of each $\Delta_i$, such that certain sums are $2 \pi$.

\begin{figure}[h]
  \begin{tikzpicture}[framed]
    \coordinate (A) at (0,0);
    \coordinate (B) at (3.8,0);
    \coordinate (C) at (4,1.8);
    \coordinate (D) at (1,4);

    \coordinate (txtstart) at ($(A)!0.5!(D)$);
    \coordinate (txtend) at ($(B)!0.5!(D)$);

    \draw (A) ++(-0.4,0) node {$0$};
    \draw (B) ++(0.4,0) node {$1$};
    \draw (C) ++(0.4,0) node {$z$};
    \draw (D) ++(-0.4,0) node {$\infty$};

    \coordinate (A1) at ($(A)!0.8!(D)$);
    \coordinate (B1) at ($(B)!0.8!(D)$);
    \coordinate (C1) at ($(C)!0.8!(D)$);
    \coordinate (A2) at ($(A)!0.8!(C)$);
    \coordinate (B2) at ($(B)!0.8!(C)$);
    \coordinate (D2) at ($(D)!0.8!(C)$);
    \coordinate (A3) at ($(A)!0.8!(B)$);
    \coordinate (C3) at ($(C)!0.8!(B)$);
    \coordinate (D3) at ($(D)!0.8!(B)$);
    \coordinate (B4) at ($(B)!0.8!(A)$);
    \coordinate (C4) at ($(C)!0.8!(A)$);
    \coordinate (D4) at ($(D)!0.8!(A)$);
    \draw (A1) -- (D4);
    \draw (A2) -- (C4);
    \draw (A3) -- (B4);
    \draw (B1) -- (D3);
    \draw (B2) -- (C3);
    \draw (C1) -- (D2);
    \filldraw[opacity=0.1] (A1) -- (B1) -- (C1) -- cycle;
    \draw (C1) -- (A1) -- (B1);
    \draw[very thick, ->] (B1) -- (C1);
    \draw[dashed] (A1) -- (D);
    \draw[dashed] (B1) -- (D);
    \draw[dashed] (C1) -- (D);
    \filldraw[opacity=0.1] (A2) -- (B2) -- (D2) -- cycle;
    \draw (A2) -- (B2) -- (D2) -- cycle;
    \draw[dashed] (A2) -- (C);
    \draw[dashed] (B2) -- (C);
    \draw[dashed] (D2) -- (C);
    \filldraw[opacity=0.1] (A3) -- (C3) -- (D3) -- cycle;
    \draw (A3) -- (C3) -- (D3) -- cycle;
    \draw[dashed] (A3) -- (B);
    \draw[dashed] (C3) -- (B);
    \draw[dashed] (D3) -- (B);
    \filldraw[opacity=0.1] (B4) -- (C4) -- (D4) -- cycle;
    \draw[very thick, ->] (B4) -- (C4);
    \draw (C4) -- (D4) -- (B4);
    \draw[dashed] (B4) -- (A);
    \draw[dashed] (C4) -- (A);
    \draw[dashed] (D4) -- (A);

    \draw[->] (txtstart) ++(-60:.4 and 0.1) arc (-60:60:.4 and 0.1);
    \draw     let \p1 = (txtstart) in let \p2 = (txtend) in
              (txtstart) ++(0.45,0) node[anchor=west,text width=\x2 - \x1,
                  scale=0.7]{$z = \crossratio{\infty}{0}{1}{z}$};

  \end{tikzpicture}
  \caption{Identifying torus traversal with multiplication in $\mathbb{C}$.}
  \label{fig:angles}
\end{figure}
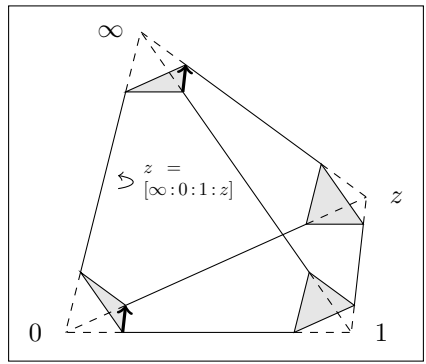

By considering each $\Delta_i$ as ideal, each segment of such a path can
be viewed as complex multiplication (recall that all cross-ratios are in
$\mathbb{C} \cup \set{\infty}$.   For example, in figure~\ref{fig:angles},
the rotation that takes a point from start to end of a bold arrow is
the same rotation that takes $1$ to $z$: multiplication by $z$. By
re-ordering the vertices of $\Delta_i$, other paths may be expressed as
re-ordered cross-ratios.

\begin{definition}
  Let a triangulation be fixed. To each simplex associate a variable
  $z_i$. For each edge $e_{ij}$ in the simplex, transforming any point
  by rotating around $e_{ij}$ is equivalent to multiplication of that
  point by either $z_i$, $\frac{z_i - 1}{z_i}$, or $\frac{1}{1-z_i}$
  (fix one angle as $z_i$, determine the rest by cross-ratio
  rearrangement).

  Let $p$ be a closed path in a triangulated surface, transverse to
  all edges. For every consecutive pair of edges that $p$ crosses,
  those edges must meet at an angle, which is associated to some
  $z_i$, $\frac{z_i - 1}{z_i}$, or $\frac{1}{1-z_i}$. The
  \definedterm{path monodromy of $p$}, denoted $\mathcal{M}(p)$, is
  the product of all these variables.
\end{definition}

\begin{definition}
  Let a triangulation of a manifold $M$ be fixed. For each torus
  boundary component $T_i$ of $M$, choose meridian and longitude paths
  in the torus (transverse to all edges) $m_i$ and $\ell_i$
  respectively. The \definedterm{cusp equations} are
  \begin{equation}
    \label{eq:cusp}
    \mathcal{M}(m_i) = 1 \qquad\text{and}\qquad \mathcal{M}(\ell_i) = 1
  \end{equation}
  There are two such equations for each torus boundary component.
\end{definition}

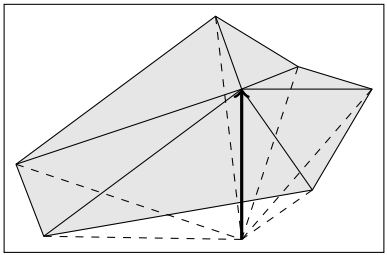
\begin{figure}[h]
  \begin{tikzpicture}[framed]
    \coordinate (base) at (0,-1);
    \coordinate (top) at (0,1);
    \coordinate (A) at (30:2);
    \coordinate (B) at (60:1.5);
    \coordinate (C) at (100:2);
    \coordinate (D) at (180:3);
    \coordinate (E) at (200:2.8);
    \coordinate (F) at (340:1);

    \filldraw[opacity=0.1] (A) -- (B) -- (C) -- (D) -- (E) -- (F) -- cycle;
    \draw (A) -- (B) -- (C) -- (D) -- (E) -- (F) -- cycle;
    \draw (A) -- (top);
    \draw (B) -- (top);
    \draw (C) -- (top);
    \draw (D) -- (top);
    \draw (E) -- (top);
    \draw (F) -- (top);
    \draw[dashed] (A) -- (base);
    \draw[dashed] (B) -- (base);
    \draw[dashed] (C) -- (base);
    \draw[dashed] (D) -- (base);
    \draw[dashed] (E) -- (base);
    \draw[dashed] (F) -- (base);
    \draw[very thick, ->] (base) -- (top);
  \end{tikzpicture}
  \caption{Piecing together a polygon from a triangulation.}
  \label{fig:sort_of_lens_space_looking_thing}
\end{figure}

One more piece of geometric information will be useful. For each edge
in the triangulation, consider all face identifications that contain
that edge as an axis. Such identifications identify the edges of
neighboring faces, and these faces join to form two polygons, one
around each vertex of the edge. For example, see figure
\ref{fig:sort_of_lens_space_looking_thing}, in which the shaded region
is a polygon for the edge shown in bold.  We arbitrarily consider the
polygon to be the one about the terminus point of the edge.

\begin{definition}
  Let a triangulation of a manifold $M$ be fixed. For each distinct
  edge $e_i$ in the triangulation, let $P_i$ be the polygon described
  above associated to that edge. Let $p_i$ be a path in $P_i$ circling
  the central point. The \definedterm{edge equations} are
  \begin{equation}
    \label{eq:edge}
    \mathcal{M}(p_i) = 1
  \end{equation}
  There are as many edge equations as there are distinct edges in the
  chosen triangulation.
\end{definition}

If two triangulations of $M$ are both composed of non-degenerate,
geometrically viable simplices, by Mostow rigidity they must yield the
same volume. Therefore, given a triangulation $\set{\Delta_i}$ of $M$,
if the cross-ratios $z_i$ all have positive imaginary part, then
computing $\sum_i D(\Delta_i)$ will yield $\vol(M)$. It is possible,
however, for a triangulation not to yield solutions that give a
geometric volume. In software implementations, this is usually dealt
with by re-triangulating using different parameters.

\begin{figure}[h]
  \begin{minipage}{.35\textwidth}
    \centering
    \includegraphics{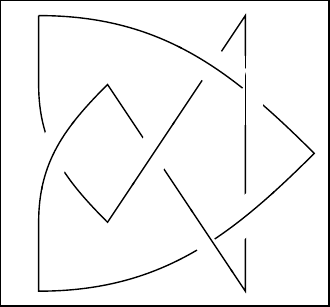}
    \captionof{figure}{$6_1$.}
    \label{fig:6_1}
  \end{minipage}
  \begin{minipage}{0.5\textwidth}
    \centering
    \includegraphics{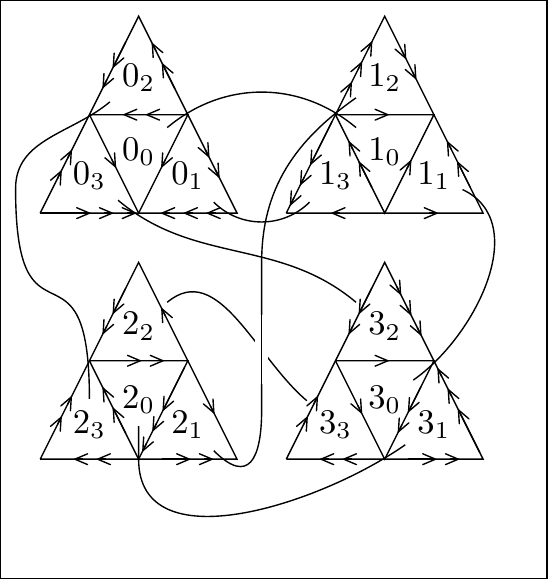}
    \captionof{figure}{$m032$.}
    \label{fig:m032_triangulation}
  \end{minipage}
\end{figure}

\begin{figure}[h]
  \includegraphics{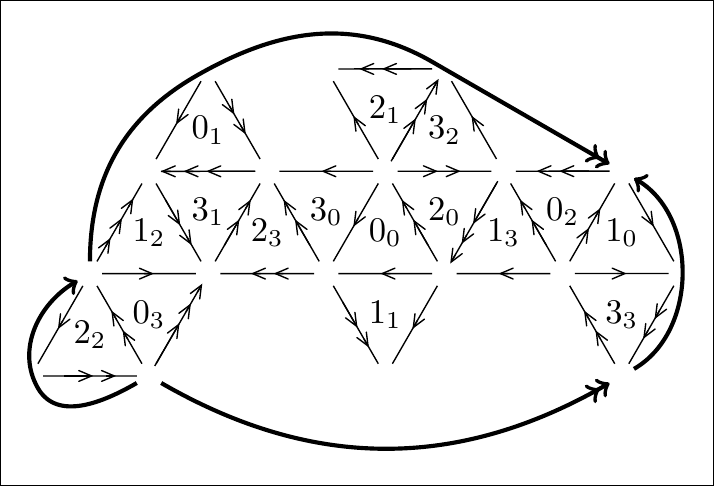}
  \caption{The torus boundary of $m032$.}
  \label{fig:m032_boundary}
\end{figure}

As presented, the system of cusp and edge equations can be hard to solve
by hand, because their degree is unbounded.  For example, the complement
of the knot $6_1$ (see figure~\ref{fig:6_1}), known as $m032$, has the
triangulation given in figure~\ref{fig:m032_triangulation}, and admits
the boundary torus given in figure~\ref{fig:m032_boundary}.

A path through the torus from left to right (for example, from the
single-arrow edge of $2_2$ to the corresponding edge of $1_0$) would
have to pass through at least eleven faces, inducing a cusp equation with
eleven terms. Computer programs such as \prog{SnapPy} \parencite{snappy}
can solve these equations (numerically) quite efficiently.

\subsection{Ptolemy coordinates}

The Ptolemy coordinates of Garoufalidis-Thurston-Zickert offer an
alternative representation of the information encoded in the gluing
and cusp equations.  For our purposes, their principal use is encoding
information about arbitrary representations
$\rho : \pi_1(M) \to \PSL(n, \mathbb{C})$.  We refer the reader
to \parencite{garoufalidis2015}, \parencite{1401.5542}, and present a
brief overview here, specialized for the particular case of hyperbolic
$3$-manifolds.

\begin{definition}
  For a fixed triangulation of a manifold $M$ (with vertices labeled
  $0$, $1$, $2$, $3$), a \definedterm{Ptolemy assignment} is an
  association of variables in $\mathbb{C}^{\times}$ to each edge of
  each simplex satisfying equations~\ref{eq:ptolemy_relations}
  below. We denote the variable of the $i$th simplex on the edge
  between vertices $j$ and $k$ as $c^i_{jk}$.

  If two edges are identified, their associated variables are the
  same.  This immediately implies, for example, that
  $c^i_{jk} = -c^i_{kj}$.

  The \definedterm{Ptolemy relations} are
  \begin{equation}
    \label{eq:ptolemy_relations}
    c_{03}^{i}c_{12}^{i} + c^i_{01}c^i_{23} = c^i_{02}c^i_{13}
  \end{equation}
\end{definition}

\begin{proposition}
  For a fixed triangulation of a manifold $M$, and for any Ptolemy
  assignment of this triangulation,
  \begin{equation}
    \label{eq:ptolemy_to_cross_ratio}
    z_i = \pm \frac{c^i_{03} c^i_{12}}{c^i_{02} c^i_{13}}
  \end{equation}
  where the sign is determined by an obstruction class, as described
  in \parencite{1401.5542} and \parencite{garoufalidis2015}.
\end{proposition}

Since in our work we are interested only in $z_i$, we may often assume
without loss of generality that certain $c^i_{jk} = 1$, where the number
of such variables depends on the number of cusps of the manifold.

In the case of $m032$ (the manifold presented above), there are only three
distinct edges, therefore three independent $c^i_{jk}$ variables. The
Ptolemy relations are homogeneous degree 2 equations in two variables,
and taking one of the variables to be $1$ makes this a readily solvable
system.

\subsection{Lattice matching}

Given a large collection of (numerically approximated) hyperbolic
3-manifold volumes (with the same invariant trace field), we desired
to construct the coarsest lattice containing the point given by each
volume. We desired to reconstruct the linear combination whose sum was
the volume of the manifold - we refer to the volume as the projection
of the lattice point to $\mathbb{R}$. Searching for representative
manifolds of the generating basis of this lattice would then provide
insight into the the relation between volumes and the Bloch group. The
question we sought to answer therefore was ``Given the dimension $d$
of a lattice, as well as a set $S$ of projections of points in that
lattice to $\mathbb{R}$, what is a $d$-subset $B$ of $S$ which, up to
a minimal scaling factor, contains all other lattices generated by
$d$-subsets of $S$?''.

\subsubsection{An assumption}
\label{sec:an-assumption}

In answering this problem, we had no exact norms: instead of $S$, we
had numerical approximations of elements of $S$. Using only numerical
approximations, we cannot prove that any two elements of $S$ are linearly
independent. However, we can be reasonably confident that two numerical
representations represent linearly independent elements.  In our work, we
used upwards of 50 places of precision and used the \prog{LLL} algorithm
implemented in \prog{PARI} \parencite{PARI2}.  If the algorithm indicated
that coefficients with magnitude greater than $2^{12}$ was required
to obtain a linear dependence, we considered the two elements linearly
independent.  As we very rarely observed coefficients with magnitudes
greater than 10 and never observed enough linearly independent elements
to contradict Borel's result below, we consider this reasonable.

\subsubsection{The lattice dimension}

By a result of Borel \parencite[397]{neumannhilberts3rd}, for an
invariant trace field $\mathbb{Q}[x]/\langle p \rangle$, where $p$ has
$r_2$ complex places, $\bloch{F}/(\text{torsion})$ is isomorphically a
lattice in $\mathbb{R}^{r_2}$. We therefore worked under the
assumption that $d = r_2$.  In some case, we observed strictly less
than $r_2$ linearly independent volumes; we believe this is due to
insufficient data, and that given unlimited resources we would
eventually discover another, linearly independent volume.

\subsubsection{An example}

Suppose $d = 2$, and we wish to obtain a $B$ for the following (numerical
approximations of) $S$:
\begin{align*}
  V &= 2.7182818284590 & (\approx e) \\
  W &= 5.4365636569180 & (\approx 2 e) \\
  X &= 6.3496623769612 & (\approx 15 \pi - 15 e) \\
  Y &= 11.7197489640976 & (\approx 2 \pi + 2 e) \\
  Z &= 21.1445269248670 & (\approx 5 \pi + 2 e).
\end{align*}

The coarsest lattice generated by $V,W,X,Y,Z$ may be generated by
$\pi$ and $e$.  We cannot obtain a $2$ element subset which generates
this lattice, but we would like to recover something close, ideally
containing $V$. We would also like to obtain a measure of how much
this $2$ element subset $B$ fails to completely generate $S$.

\begin{definition}
  For a $d$-subset $B$ of $S$, as described above, the \definedterm{fit
  ratio} of $B$ is the smallest positive integer $f$ such that the
  lattice with generators $B$ contains $fS$. If no such integer exists,
  the fit ratio of $B$ is $\infty$.
\end{definition}

We certainly do not wish to select the two smallest elements of $S$:
in the example $\set{V,W}$ has fit ratio $\infty$.  Nor is it
sufficient to simply select the two least linearly independent
elements: $\set{V, X}$ has fit ratio $15$ by considering $Y$.

A good choice would be $B = \set{V, Y}$, with fit ratio
$2$.  Our algorithm for detecting such $B$ is formalized in
algorithm~\ref{alg:lattice-match}, which is implemented in our
experimentation software.

\subsubsection{An algorithm}

To find a good choice of $B$, we first find some linearly independent
basis of $S$ (by brute force, if necessary), then express each element
of $S$ in terms of that basis, yielding vectors in $\mathbb{Q}^d$
(if this could not be achieved, it would be a contradiction of Borel's
result). We then test all $d$-sized sets of vectors, and select the one
with the least determinant. The ratio of this determinant to the
fractional GCD of all such determinants yields the fit ratio.

\begin{algorithm}
  \caption{Best-fit lattice matching algorithm.}
  \label{alg:lattice-match}
  \begin{algorithmic}
    \State $manifolds \gets \Call{sort\_by\_increasing\_volume}{manifolds}$
    \State $next\_i \gets 2$ \Comment{Which $e_i$ we seek a manifold for}
    \State $rational\_vecs \gets \set{(manifolds_1, e_1)}$
        \Comment{Tuples ($M$, $v$), with $v \in \mathbb{Q}^n$}
    \State $rational\_basis \gets \set{manifolds_1}$
        \Comment{All $M$ associated to some $e_i$}
    \ForAll{$M \in manifolds$ except $manifolds_1$}
      \State $dep = \Call{lindep}{\vol(M),rational\_basis}$
          \Comment{Via PARI; $\langle dep, volumes\rangle = 0$}
      \If{$dep$ shows $\vol(M)$ is not a linear combination of others}
        \If{$dep$ shows $\vol(M)$ is a rational combination of others}
           \State $rational\_vecs \gets rational\_vecs \cup \set{(M), dep_1^{-1}
                   dep_{2,3,\dotsc,n+1}}$
        \Else
           \If{$next\_i > n$}
             \State \Return \text{Failure}
                 \Comment{There are more than $n$ l.i.\ manifolds}
           \EndIf
           \State $rational\_basis \gets rational\_basis \cup \set{M}$
           \State $rational\_vecs \gets rational\_vecs \cup \set{(M),e_{next\_i}}$
           \State $next\_i \gets next\_i + 1$
        \EndIf
      \EndIf
    \EndFor
    \State $gcd \gets \infty$
    \State $best\_det \gets \infty$
    \State $best\_basis \gets \varnothing$
    \ForAll{$B \subseteq rational\_vecs$ with $|B| = n$}
      \State $det \gets \det(\set{v \given (M,v) \in B})$
      \If{$det \ne 0$}
        \State $gcd \gets \operatorname{fractional\_gcd}(gcd, det)$
        \If{$det < best\_det$}
          \State $best\_det \gets det$
          \State $best\_basis \gets B$
        \EndIf
      \EndIf
    \EndFor
    \State \Return (basis, fit ratio) $:=(best\_basis, \frac{best\_det}{gcd})$
  \end{algorithmic}
\end{algorithm}

In our implementation, we use \prog{PARI} for the \codebit{lindep}
algorithm, and as noted in \ref{sec:an-assumption}, we use a combination
of high precision and detection of large coefficients to determine
linear dependence.

Applied to the example, \codebit{rational\_vecs} is calculated as \[
\set*{(V, (1,0)), (X, (0,1)), (Y, (\frac{1}{4}, -\frac{2}{15})), (Z,
(-7, -\frac{1}{3}))}, \] discarding $W$ since it is a linear multiple
of $V$. Calculating determininants of potential bases gives:
\begin{align*}
  \set{V,X} \rightarrow \det \begin{bmatrix}1 & 0 \\ 0 & 1\end{bmatrix} &= 1 & \set{V,Y} \rightarrow \det \begin{bmatrix}1 & 4 \\ 0 & -\frac{2}{15}\end{bmatrix} &= -\frac{2}{15} \\
  \set{V,Z} \rightarrow \det \begin{bmatrix}1 & -7 \\ 0 & -\frac{1}{3}\end{bmatrix} &= -\frac{1}{3} & \set{X,Y} \rightarrow \det \begin{bmatrix}0 & 4 \\ 1 & -\frac{2}{15}\end{bmatrix} &= 4 \\
  \set{X,Z} \rightarrow \det \begin{bmatrix}0 & -7 \\ 1 & -\frac{1}{3}\end{bmatrix} &= 7 & \set{Y,Z} \rightarrow \det \begin{bmatrix}4 & -7 \\ -\frac{2}{15} & -\frac{1}{3}\end{bmatrix} &= -\frac{34}{15}.
\end{align*}
        
Since $\frac{2}{15}$ is the smallest absolute value of all possible
determinants, $B = \set{V, Y}$ is returned by the algorithm. Since the
fractional GCD of all possible determinants is $\frac{1}{15}$, the fit
ratio of $B$ is $2$.

\section{Experimentation}
\label{sec:main-result}

Our main results are the following data.

For our work, we used the volumes for a large number of hyperbolic
3-manifolds, which we obtained by performing Dehn surgery on manifolds
from a base census.  We used \prog{SnapPy} \parencite{snappy}, although
for invariant trace field calculations we used the older \prog{Snap}
\parencite{snap} interface.

\subsection{Manifold generation}

The code we used to collect manifold information can be found at
\url{https://github.com/s-gilles/maps-reu-code}. We considered all
orientable, cusped, hyperbolic manifolds which can be triangulated
using 9 tetrahedra or less (the \codebit{OrientableCuspedCensus} within
\prog{SnapPy} at the time of work), as well as all link complements
using 3 crossings or more (a subset of \codebit{LinkExteriors}). We
also considered additional manifolds in the \codebit{LinkExteriors}
and \codebit{HTLinkExteriors} collections.

Of these manifolds, we performed Dehn surgery of type $(p,q)$ over
each of the $n$ cusps, with $p \in [0, L(n)]$ and
$q \in [-L(n), L(n)]$, (and $(p,q) = 1$) where $L(n)$ is given by
figure~\ref{fig:ln}.

\begin{figure}[h]
  \begin{tabular}{|l|cccccccc|} \hline $L(n)$ & $16$ & $12$ & $8$ &
    $6$ & $4$ & $3$ & $3$ & $2$ \\ \hline $n $ & $ 1$ & $ 2$ & $3$ &
    $4$ & $5$ & $6$ & $7$ & $ \ge 8$ \\ \hline
  \end{tabular}
  \caption{$L(n)$}
  \label{fig:ln}
\end{figure}

Of the resulting manifolds, we discarded those whose invariant trace
field resulted in a polynomial with degree $9$ or more (or if
\prog{SnapPy} or \prog{Snap} were unable to triangulate the
result). For each remaining manifold, we we stored intermediate
information about the manifold.  The resulting file, available at
\url{http://www.curve.unhyperbolic.org/linComb/finalized_data/volumes/all_volumes.csv},
contains over 790,000 manifolds, representing over 6,300 distinct
invariant trace fields, with volumes accurate to at least 50 decimal
places. The work was performed on the University of Maryland's
computation cluster.

\subsection{Lattice matching}

With the data stored in \codebit{all\_volumes.csv}, we have enough
data to conjecture partial lattice generators for 5,900 invariant
trace fields. Of those, for 312 invariant trace fields we have
enough data to conjecture not only generating volumes for the full
lattice, but also representative manifolds associated to the invariant
trace field which exhibit those volumes. This data can be found at
\url{http://www.curve.unhyperbolic.org/linComb/finalized_data/volume_spans.csv},
and the code we used for conjecturing lattice generators is available
at \url{https://github.com/s-gilles/maps-reu-code}, an implementation
of algorithm~\ref{alg:lattice-match}.

Sample data can be found in section~\ref{sec:lattice-examples}.

\subsection{Linear combinations for exotic volumes} While we used the
Ptolemy coordinates for computing volumes of hyperbolic 3-manifolds,
the computations may be generalized to compute the volume of any
representation $\rho: \Gamma \to \PSL(n, \mathbb{C})$.

The extended Bloch group is exactly the Bloch group up
to torsion, so another method for experimentally testing
conjecture~\ref{conj:neumann-main} is to produce linear combinations of
hyperbolic volumes matching exotic volumes, where these exotic volumes
arise from generalized representations into $\PSL(n,\mathbb{C})$ where
$n$ is not necessarily $2$.

Using the data of \codebit{all\_volumes.csv}, we have been able
to find a great many such combinations.  This data is available at
\url{http://www.curve.unhyperbolic.org/linComb/finalized_data/linear_combinations/}.

Sample data can be found in section~\ref{sec:lattice-examples}.

\section{A counter-example to
  conjecture~\ref{conj:false-strengthening}}

\label{sec:strengthening-counterexample} A secondary result of our work
is a counter-example to conjecture~\ref{conj:false-strengthening}.
Assuming the conjecture holds, two of the manfold volumes catologued
would imply the existence of a manifold with an impossibly small volume.

\begin{theorem}[{\parencite[Corollary~1.3]{gmm-min-vol}}]
  \label{thm:weeks-is-min}
  The Weeks manifold is the unique closed orientable hyperbolic 3-manifold
  of smallest volume: $0.9427\dotsc$.
\end{theorem}

\begin{counterexample}
  Conjecture~\ref{conj:false-strengthening} is false.
\end{counterexample}
\begin{proof}
  The polynomial $p(x) = x^3 - x^2 + 1$ has exactly one complex
  place. The manifold $m003(2,1)$ has hyperbolic volume
  $0.9427\dotsc$, and the manifold $10^2_{94}(2,3)(5,11)$ has
  hyperbolic volume $1.4140\dotsc$.

  If conjecture~\ref{conj:false-strengthening} were true, the volumes of
  these manifolds (which have invariant trace field $\mathbb{Q}[x]/\langle
  p \rangle$) would be expressable as a linear combination of one
  generator: $\vol(M)$. Then \[ \vol(M) \mid \frac{0.9427\dotsc}{2} =
  (vol(10^2_{94}(2,3)(5,11)) - \vol(m003(2,1))). \] This would contradict
  theorem~\ref{thm:weeks-is-min}.
\end{proof}

In our entire census, we did not observe any counterexamples where the
obstruction factor was less than $\frac{1}{2}$.  We suspect this is
related to the following theorem of Neumann.

\begin{theorem}[{\parencite[Theorem~2.7]{neumann06}}]
  If $M$ has cusps then $[M]$ is defined in $\bloch{F}$, for $F$ the
  invariant trace field of $M$, while if $M$ is closed $2[M]$ is
  defined in $\bloch{F}$.
\end{theorem}

\section{Evidence for conjecture~\ref{conj:neumann-weaker}}
\label{sec:manifold-existence}

Our data provides an opportunity for examining the strength of
conjecture~\ref{conj:neumann-weaker}. It is known that every field with
one complex place arises as the invariant trace field of a hyperbolic
manifold \parencite{maclachlanreid}, but an open question in general
\parencite{neumann06}.

\subsection{Examples of missing fields}

Our data is necessarily incomplete.  Some interesting absences:

\begin{enumerate}

\item Some fields associated to ``simple'' polynomials such as $x^4 -
2x^2 + 4$, $x^4 + 9$, and $x^7 - 3$ do not arise.  This also includes
some degree $2$ polynomials such as $x^2 + 13$, $x^2 + 22$, and $x^2 -
x + 17$, even though it is known that all fields with one complex place
must be exhibited by some manifold.

\item Some number fields, such as the two concrete number fields
arising from $\mathbb{Q}[x]/\langle x^4 - x + 2 \rangle$, produced
surprisingly little data.

For the concrete number field associated to root $-0.850\dotsc \pm
1.01\dotsc i$, we observed only two volumes: $21.531\dotsc$ and
$16.383\dotsc$. These two volumes are linearly independent, but we know
very little about the proposed lattice.

For the concrete number field associated to root $0.850\dotsc
\pm 0.654\dotsc i$, we only observed one volume: $9.054\dotsc$,
though the lattice should be of dimension $2$. This is (from some
perspective) the simplest example from our data which does not support
conjecture~\ref{conj:neumann-main}.

\end{enumerate}

\subsection{Data analysis}

By comparing the fields we observed (under some reasonable restrictions) to a
proven-complete census under the same restrictions, we obtain an estimate
of the probability that an arbitrary field (within these restrictions)
appears in our census. If conjecture~\ref{conj:neumann-weaker} holds,
then our census samples from all possible fields. If our sample was
uniformly random, we would then expect that our census contains the same
percentage of fields with $r_2 = 1$ as for $r_2 \ne 1$.

For our calculations, we used a online census \parencite{jones2014}
of abstract number fields expressed as polynomials. For
the field restrictions we used, this census has been proven
complete. In figures~\ref{fig:observed-concrete-fields} and
\ref{fig:observed-abstract-fields}, $n$ refers to the degree of
the polynomial $p$ to which the field is associated, and $D$ is the
discriminant of the polynomial. We only consider restrictions for which
multiple values of $r_2$ are possible, and we ignore $r_2 = 4$ due to
relative scarcity of data (the total fields for $r_2 = 4$ are so numerous
that our results are scarce enough to be statistical noise).

\begin{figure}[h]
  \centering
  {\small
  \begin{tabular}{|ll|lll|lll|lll|} \hline
    \multicolumn{2}{|l|}{Restrictions} & \multicolumn{3}{c|}{$r_2 = 1$} & \multicolumn{3}{c|}{$r_2 = 2$} & \multicolumn{3}{c|}{$r_2 = 3$} \\
    $n$ & $|D|^{1/n}$ & {\small Found} & {\small Total} & \% & {\small Found} & {\small Total} & \% & {\small Found} & {\small Total} & \% \\ \hline
    $4$ & $\le 8$ & $56$ & $137$ & $40.8\%$ & $76$ & $408$ & $18.6\%$ & & & \\
    $4$ & $\le 10$ & $65$ & $444$ & $14.6\%$ & $82$ & $1100$ & $7.4\%$ & & & \\
    $4$ & $\le 12$ & $66$ & $1056$ & $6.2\%$ & $84$ & $2550$ & $3.2\%$ & & & \\
    $4$ & $\le 15$ & $66$ & $3069$ & $2.1\%$ & $84$ & $6728$ & $1.2\%$ & & & \\
    $5$ & $\le 8$ & $47$ & $77$ & $61.0\%$ & $226$ & $736$ & $30.7\%$ & & & \\
    $5$ & $\le 10$ & $65$ & $472$ & $13.7\%$ & $326$ & $3470$ & $9.3\%$ & & & \\
    $5$ & $\le 12$ & $73$ & $1670$ & $4.3\%$ & $356$ & $10992$ & $3.2\%$ & & & \\
    $5$ & $\le 15$ & $76$ & $7556$ & $1.0\%$ & $362$ & $41776$ & $.8\%$ & & & \\
    $6$ & $\le 8$ & $11$ & $40$ & $27.5\%$ & $180$ & $1222$ & $14.7\%$ & $166$ & $1851$ & $8.9\%$ \\
    $6$ & $\le 10$ & $28$ & $405$ & $6.9\%$ & $374$ & $8434$ & $4.4\%$ & $291$ & $10887$ & $2.6\%$ \\
    $6$ & $\le 12$ & $37$ & $2335$ & $1.5\%$ & $514$ & $38722$ & $1.3\%$ & $352$ & $42123$ & $.8\%$ \\
    $6$ & $\le 15$ & $48$ & $15556$ & $.3\%$ & $578$ & $204108$ & $.2\%$ & $380$ & $190995$ & $.1\%$ \\
    $7$ & $\le 10$ & $10$ & $137$ & $7.2\%$ & $276$ & $8070$ & $3.4\%$ & $391$ & $38103$ & $1.0\%$ \\
    $7$ & $\le 12$ & $14$ & $1473$ & $.9\%$ & $455$ & $57292$ & $.7\%$ & $637$ & $219879$ & $.2\%$ \\
    $7$ & $\le 15$ & $30$ & $16759$ & $.1\%$ & $599$ & $506188$ & $.1\%$ & $890$ & $1612152$ & $.0\%$ \\
    $8$ & $\le 10$ & $1$ & $22$ & $4.5\%$ & $75$ & $752$ & $9.9\%$ & $288$ & $3141$ & $9.1\%$ \\
    $8$ & $\le 12$ & $7$ & $246$ & $2.8\%$ & $199$ & $5808$ & $3.4\%$ & $584$ & $16764$ & $3.4\%$ \\
    $8$ & $\le 15$ & $11$ & $2560$ & $.4\%$ & $383$ & $50268$ & $.7\%$ & $976$ & $120111$ & $.8\%$ \\
    \hline
  \end{tabular}
  }
  \caption{Observed concrete field percentages by restriction}
  \label{fig:observed-concrete-fields}
\end{figure}

\begin{figure}[h]
  \centering
  {\small
  \begin{tabular}{|ll|lll|lll|lll|} \hline
    \multicolumn{2}{|l|}{Restrictions} & \multicolumn{3}{c|}{$r_2 = 1$} & \multicolumn{3}{c|}{$r_2 = 2$} & \multicolumn{3}{c|}{$r_2 = 3$} \\
    $n$ & $|D|^{1/n}$ & {\small Found} & {\small Total}& \% & {\small Found} & {\small Total}& \% & {\small Found} & {\small Total}& \% \\ \hline
    $4$ & $\le 8$ & $56$ & $137$ & $40.8\%$ & $60$ & $204$ & $29.4\%$ & & & \\
    $4$ & $\le 10$ & $65$ & $444$ & $14.6\%$ & $66$ & $550$ & $12.0\%$ & & & \\
    $4$ & $\le 12$ & $66$ & $1056$ & $6.2\%$ & $68$ & $1275$ & $5.3\%$ & & & \\
    $4$ & $\le 15$ & $66$ & $3069$ & $2.1\%$ & $68$ & $3364$ & $2.0\%$ & & & \\
    $5$ & $\le 8$ & $47$ & $77$ & $61.0\%$ & $168$ & $368$ & $45.6\%$ & & & \\
    $5$ & $\le 10$ & $65$ & $472$ & $13.7\%$ & $264$ & $1735$ & $15.2\%$ & & & \\
    $5$ & $\le 12$ & $73$ & $1670$ & $4.3\%$ & $294$ & $5496$ & $5.3\%$ & & & \\
    $5$ & $\le 15$ & $76$ & $7556$ & $1.0\%$ & $300$ & $20888$ & $1.4\%$ & & & \\
    $6$ & $\le 8$ & $11$ & $40$ & $27.5\%$ & $162$ & $611$ & $26.5\%$ & $143$ & $617$ & $23.1\%$ \\
    $6$ & $\le 10$ & $28$ & $405$ & $6.9\%$ & $350$ & $4217$ & $8.2\%$ & $263$ & $3629$ & $7.2\%$ \\
    $6$ & $\le 12$ & $37$ & $2335$ & $1.5\%$ & $489$ & $19361$ & $2.5\%$ & $324$ & $14041$ & $2.3\%$ \\
    $6$ & $\le 15$ & $48$ & $15556$ & $.3\%$ & $553$ & $102054$ & $.5\%$ & $352$ & $63665$ & $.5\%$ \\
    $7$ & $\le 10$ & $10$ & $137$ & $7.2\%$ & $271$ & $4035$ & $6.7\%$ & $383$ & $12701$ & $3.0\%$ \\
    $7$ & $\le 12$ & $14$ & $1473$ & $.9\%$ & $450$ & $28646$ & $1.5\%$ & $629$ & $73293$ & $.8\%$ \\
    $7$ & $\le 15$ & $30$ & $16759$ & $.1\%$ & $594$ & $253094$ & $.2\%$ & $882$ & $537384$ & $.1\%$ \\
    $8$ & $\le 10$ & $1$ & $22$ & $4.5\%$ & $75$ & $376$ & $19.9\%$ & $286$ & $1047$ & $27.3\%$ \\
    $8$ & $\le 12$ & $7$ & $246$ & $2.8\%$ & $199$ & $2904$ & $6.8\%$ & $582$ & $5588$ & $10.4\%$ \\
    $8$ & $\le 15$ & $11$ & $2560$ & $.4\%$ & $383$ & $25134$ & $1.5\%$ & $974$ & $40037$ & $2.4\%$ \\
    \hline
  \end{tabular}
  }
  \caption{Observed abstract field percentages by restriction}
  \label{fig:observed-abstract-fields}
\end{figure}

Figure~\ref{fig:observed-concrete-fields} implies, in general, a
negative correlation between the number of complex places of the
field and the probability that it is exhibited by some manifold.
Figure~\ref{fig:observed-abstract-fields} offers hope for an explanation.

While figure~\ref{fig:observed-concrete-fields}
lists the observed percentage of concrete number fields,
figure~\ref{fig:observed-abstract-fields} lists the observed percentage
of fields considered as $\mathbb{Q}[x]/\langle p \rangle$: in particular,
we treat any two concrete fields arising from different roots of the
same polynomial as identical. The percentages are relatively even
across rows in this case, which supports that our census was selecting
from all possible abstract number fields. This partially supports
conjecture~\ref{conj:neumann-weaker}, and more solidly supports the
following weaker version.

\begin{conjecture}
  \label{conj:only-abstract-represented}
  Every non-real abstract number field $\mathbb{Q}[x]/\langle p \rangle$
  arises as abstract number field associated to the invariant trace
  field of some hyperbolic manifold.
\end{conjecture}

\section{Selected examples}
\label{sec:lattice-examples}

We present a selection of lattices, together with the manifolds which may
be generators.  These are constructed from the \codebit{volume\_spans.csv}
and \codebit{linear\_combinations} files of our results.

Each of the following lattices are two-dimensional: each is associated
to a polynomial with two complex places. Our results contain data for
three and even four complex places, but these are harder to represent
visually. In the following graphs, one potential generator is arbitrarily
assigned to the $x$-axis, and the other to the $y$-axis.  A point at
position $(x,y)$ represents a volume $v = a \vol(g_1) + b \vol(g_2)$,
where $x = a \vol(g_1)$ and $y = b \vol g_2)$.

The gray dots indicate our predicted lattice (based on the volumes
of all manifolds associated to that invariant trace field), while the
crosses indicate manifolds we have observed: blue $+$ markers indicate
geometric volumes, while red $\times$ markers indicate exotic volumes.
The relevant volume (whether geometric or exotic) of the manifold is
the projection to $\mathbb{R}$ of the point in the lattice.

\begin{figure}[h]
  \begin{tikzpicture}[framed]
    \begin{axis}[%
                 xlabel={Projection to generator $\vol(v1859(-1,3))$},%
                 ylabel={Projection to generator $\vol(t07828)$},%
                ]
      \addplot[color=gray, mark=o, mark size=0.4pt, only marks] table {
        x       y      
        -40.48  -31.65 
        -40.48  -25.32 
        -40.48  -18.99 
        -40.48  -12.66 
        -40.48  -6.330 
        -40.48  0.0    
        -40.48  6.3306 
        -40.48  12.661 
        -40.48  18.991 
        -40.48  25.322 
        -40.48  31.653 
        -40.48  37.983 
        -40.48  44.314 
        -35.98  -31.65 
        -35.98  -25.32 
        -35.98  -18.99 
        -35.98  -12.66 
        -35.98  -6.330 
        -35.98  0.0    
        -35.98  6.3306 
        -35.98  12.661 
        -35.98  18.991 
        -35.98  25.322 
        -35.98  31.653 
        -35.98  44.314 
        -31.49  -31.65 
        -31.49  -25.32 
        -31.49  -18.99 
        -31.49  -12.66 
        -31.49  -6.330 
        -31.49  0.0    
        -31.49  6.3306 
        -31.49  12.661 
        -31.49  18.991 
        -31.49  25.322 
        -31.49  31.653 
        -31.49  37.983 
        -31.49  44.314 
        -26.99  -31.65 
        -26.99  -25.32 
        -26.99  -18.99 
        -26.99  -12.66 
        -26.99  -6.330 
        -26.99  0.0    
        -26.99  6.3306 
        -26.99  12.661 
        -26.99  18.991 
        -26.99  25.322 
        -26.99  31.653 
        -26.99  37.983 
        -26.99  44.314 
        -22.49  -31.65 
        -22.49  -25.32 
        -22.49  -18.99 
        -22.49  -12.66 
        -22.49  -6.330 
        -22.49  0.0    
        -22.49  6.3306 
        -22.49  12.661 
        -22.49  18.991 
        -22.49  25.322 
        -22.49  31.653 
        -22.49  37.983 
        -22.49  44.314 
        -17.99  -31.65 
        -17.99  -25.32 
        -17.99  -18.99 
        -17.99  -12.66 
        -17.99  -6.330 
        -17.99  0.0    
        -17.99  6.3306 
        -17.99  12.661 
        -17.99  25.322 
        -17.99  31.653 
        -17.99  37.983 
        -17.99  44.314 
        -13.49  -31.65 
        -13.49  -25.32 
        -13.49  -18.99 
        -13.49  -12.66 
        -13.49  -6.330 
        -13.49  0.0    
        -13.49  6.3306 
        -13.49  12.661 
        -13.49  18.991 
        -13.49  25.322 
        -13.49  31.653 
        -13.49  37.983 
        -13.49  44.314 
        -8.997  -31.65 
        -8.997  -25.32 
        -8.997  -18.99 
        -8.997  -12.66 
        -8.997  -6.330 
        -8.997  0.0    
        -8.997  6.3306 
        -8.997  12.661 
        -8.997  18.991 
        -8.997  25.322 
        -8.997  31.653 
        -8.997  37.983 
        -8.997  44.314 
        -4.498  -31.65 
        -4.498  -25.32 
        -4.498  -18.99 
        -4.498  -12.66 
        -4.498  -6.330 
        -4.498  0.0    
        -4.498  6.3306 
        -4.498  12.661 
        -4.498  18.991 
        -4.498  25.322 
        -4.498  31.653 
        -4.498  37.983 
        -4.498  44.314 
        0.0     -31.65 
        0.0     -25.32 
        0.0     -18.99 
        0.0     -12.66 
        0.0     -6.330 
        0.0     0.0    
        0.0     18.991 
        0.0     25.322 
        0.0     31.653 
        0.0     37.983 
        0.0     44.314 
        4.4986  -31.65 
        4.4986  -25.32 
        4.4986  -18.99 
        4.4986  -12.66 
        4.4986  -6.330 
        4.4986  6.3306 
        4.4986  12.661 
        4.4986  18.991 
        4.4986  25.322 
        4.4986  31.653 
        4.4986  37.983 
        4.4986  44.314 
        8.9973  -31.65 
        8.9973  -25.32 
        8.9973  -18.99 
        8.9973  -12.66 
        8.9973  6.3306 
        8.9973  12.661 
        8.9973  18.991 
        8.9973  25.322 
        8.9973  31.653 
        8.9973  37.983 
        8.9973  44.314 
        13.496  -31.65 
        13.496  -25.32 
        13.496  -18.99 
        13.496  -6.330 
        13.496  0.0    
        13.496  6.3306 
        13.496  12.661 
        13.496  18.991 
        13.496  25.322 
        13.496  31.653 
        13.496  37.983 
        13.496  44.314 
        17.994  -31.65 
        17.994  -25.32 
        17.994  -18.99 
        17.994  -6.330 
        17.994  0.0    
        17.994  6.3306 
        17.994  12.661 
        17.994  18.991 
        17.994  25.322 
        17.994  31.653 
        17.994  37.983 
        17.994  44.314 
        22.493  -31.65 
        22.493  -25.32 
        22.493  -18.99 
        22.493  -12.66 
        22.493  -6.330 
        22.493  0.0    
        22.493  6.3306 
        22.493  12.661 
        22.493  18.991 
        22.493  25.322 
        22.493  31.653 
        22.493  37.983 
        22.493  44.314 
        26.992  -31.65 
        26.992  -12.66 
        26.992  -6.330 
        26.992  0.0    
        26.992  6.3306 
        26.992  12.661 
        26.992  18.991 
        26.992  25.322 
        26.992  31.653 
        26.992  37.983 
        26.992  44.314 
        31.490  -31.65 
        31.490  -25.32 
        31.490  -18.99 
        31.490  -12.66 
        31.490  -6.330 
        31.490  0.0    
        31.490  6.3306 
        31.490  12.661 
        31.490  18.991 
        31.490  25.322 
        31.490  31.653 
        31.490  37.983 
        31.490  44.314 
      };
      \addplot[color=blue, mark=+, only marks, %
               point meta=explicit symbolic, nodes near coords,%
               visualization depends on=\thisrow{angle} \as \angle,%
               every node near coord/.style={font=\tiny, anchor=\angle}]%
        table [meta index=2] {
        x       y       label                           angle
        8.9973  0.0     $L7a1$                          200
        0.0     6.3306  $t07828$                        270
        0.0     12.661  $10^{3}_{27}$                   270
        4.4986  0.0     $v1859(-1,3)$                   15
      };
      \addplot[color=red, mark=x, only marks, %
               point meta=explicit symbolic, nodes near coords,%
               visualization depends on=\thisrow{angle} \as \angle,%
               every node near coord/.style={font=\tiny, anchor=\angle}]%
        table [meta index=2] {
        x       y       label                           angle
        17.994  -12.66  $s781$                          225
        26.992  -25.32  $L7a1$                          270
        -35.98  37.983  $10^{3}_{27}$                   270
        13.496  -12.66  $t08046$                        15
        26.992  -18.99  $o9_{39892}$                    270
        8.9973  -6.330  $m306$                          200
        -17.99  18.991  $t07828$                        270
      };
    \end{axis}
  \end{tikzpicture}

  ~

  \begin{tabular}{llll} \hline
    Manifold & Volume & Type & As linear combination \\ \hline
    $v1859(-1,3)$ & $4.4986\dotsc$ & geometric & $1\cdot\vol(v1859(-1,3))+0\cdot\vol(t07828)$ \\
    $t07828$ & $6.3306\dotsc$ & geometric & $0\cdot\vol(v1859(-1,3))+1\cdot\vol(t07828)$ \\
    $L7a1$ & $8.9973\dotsc$ & geometric & $2\cdot\vol(v1859(-1,3))+0\cdot\vol(t07828)$ \\
    $10^{3}_{27}$ & $12.661\dotsc$ & geometric & $0\cdot\vol(v1859(-1,3))+2\cdot\vol(t07828)$ \\
    $m306$ & $2.6667\dotsc$ & exotic & $2\cdot\vol(v1859(-1,3))-1\cdot\vol(t07828)$ \\
    $s781$ & $5.3334\dotsc$ & exotic & $4\cdot\vol(v1859(-1,3))-2\cdot\vol(t07828)$ \\
    $o9_{39892}$ & $8.0002\dotsc$ & exotic & $6\cdot\vol(v1859(-1,3))-3\cdot\vol(t07828)$ \\
    $10^{3}_{27}$ & $1.9942\dotsc$ & exotic & $-8\cdot\vol(v1859(-1,3))+6\cdot\vol(t07828)$ \\
    $L7a1$ & $1.6696\dotsc$ & exotic & $6\cdot\vol(v1859(-1,3))-4\cdot\vol(t07828)$ \\
    $t07828$ & $0.9971\dotsc$ & exotic & $-4\cdot\vol(v1859(-1,3))+3\cdot\vol(t07828)$ \\
    $t08046$ & $0.8348\dotsc$ & exotic & $3\cdot\vol(v1859(-1,3))-2\cdot\vol(t07828)$ \\
  \end{tabular}
  \caption{$p(x) = x^4-3x^2+4$, root $-1.322\dotsc+0.5i$, prospective
  basis: $\vol(v1859(-1,3))$, $\vol(t07828)$ with fit ratio 1.}
  \label{fig:v1859-lattice}
\end{figure}

Figure~\ref{fig:v1859-lattice} represents a typical ``good'' lattice
from our data. We obtained a substantial number of distinct volumes
(each volume is recorded many times in our census), and we have a number
of pairs of volumes which differ by exactly one of our conjectured
basis elements.

\begin{figure}[h]
  \begin{tikzpicture}[framed]
    \begin{axis}[%
                 xlabel={Projection to generator $\vol(v3318)$},%
                 ylabel={Projection to generator $\vol(L14n54610)$},%
                ]
      \addplot[color=gray, mark=o, mark size=0.4pt, only marks] table {
        x       y      
        -51.60  -47.06 
        -51.60  -31.37 
        -51.60  -15.68 
        -51.60  0.0    
        -51.60  15.688 
        -51.60  31.376 
        -51.60  47.064 
        -51.60  62.752 
        -45.15  -47.06 
        -45.15  -31.37 
        -45.15  -15.68 
        -45.15  0.0    
        -45.15  15.688 
        -45.15  31.376 
        -45.15  62.752 
        -38.70  -47.06 
        -38.70  -31.37 
        -38.70  -15.68 
        -38.70  0.0    
        -38.70  15.688 
        -38.70  31.376 
        -38.70  47.064 
        -38.70  62.752 
        -32.25  -47.06 
        -32.25  -31.37 
        -32.25  -15.68 
        -32.25  0.0    
        -32.25  15.688 
        -32.25  31.376 
        -32.25  47.064 
        -32.25  62.752 
        -25.80  -47.06 
        -25.80  -31.37 
        -25.80  -15.68 
        -25.80  0.0    
        -25.80  15.688 
        -25.80  47.064 
        -25.80  62.752 
        -19.35  -47.06 
        -19.35  -31.37 
        -19.35  -15.68 
        -19.35  0.0    
        -19.35  15.688 
        -19.35  31.376 
        -19.35  47.064 
        -19.35  62.752 
        -12.90  -47.06 
        -12.90  -31.37 
        -12.90  -15.68 
        -12.90  0.0    
        -12.90  31.376 
        -12.90  47.064 
        -12.90  62.752 
        -6.450  -47.06 
        -6.450  -31.37 
        -6.450  -15.68 
        -6.450  0.0    
        -6.450  31.376 
        -6.450  47.064 
        -6.450  62.752 
        0.0     -47.06 
        0.0     -31.37 
        0.0     -15.68 
        0.0     0.0    
        0.0     15.688 
        0.0     31.376 
        0.0     47.064 
        0.0     62.752 
        6.4506  -47.06 
        6.4506  -31.37 
        6.4506  -15.68 
        6.4506  15.688 
        6.4506  31.376 
        6.4506  47.064 
        6.4506  62.752 
        12.901  -47.06 
        12.901  -31.37 
        12.901  -15.68 
        12.901  0.0    
        12.901  15.688 
        12.901  31.376 
        12.901  47.064 
        12.901  62.752 
        19.352  -47.06 
        19.352  -31.37 
        19.352  0.0    
        19.352  15.688 
        19.352  31.376 
        19.352  47.064 
        19.352  62.752 
        25.802  -47.06 
        25.802  -31.37 
        25.802  -15.68 
        25.802  0.0    
        25.802  15.688 
        25.802  31.376 
        25.802  47.064 
        25.802  62.752 
        32.253  -47.06 
        32.253  -15.68 
        32.253  0.0    
        32.253  15.688 
        32.253  31.376 
        32.253  47.064 
        32.253  62.752 
        38.704  -47.06 
        38.704  -31.37 
        38.704  -15.68 
        38.704  0.0    
        38.704  15.688 
        38.704  31.376 
        38.704  47.064 
        38.704  62.752 
      };
      \addplot[color=blue, mark=+, only marks, %
               point meta=explicit symbolic, nodes near coords,%
               visualization depends on=\thisrow{angle} \as \angle,%
               every node near coord/.style={font=\tiny, anchor=\angle}]%
        table [meta index=2] {
        x       y       label                           angle
        6.4506  0.0     $v3318$                         270
        0.0     15.69   $L14n534610$                    165
      };
      \addplot[color=red, mark=x, only marks, %
               point meta=explicit symbolic, nodes near coords,%
               visualization depends on=\thisrow{angle} \as \angle,%
               every node near coord/.style={font=\tiny, anchor=\angle}]%
        table [meta index=2] {
        x       y       label                           angle
        19.352  -15.68  $K12n809$                       15
        25.802  -23.53  $t09825$                        15
        -6.450  15.688  $L13n5993$                      270
        -25.80  31.376  $10^{3}_{5}$                    270
        32.253  -31.37  $v3318$                         15
        12.901  -7.844  $t09825$                        195
        -12.90  15.688  $v3548$                         15
        -45.15  47.064  $L13n5993$                      270
        9.6760  -7.844  $v2489$                         15
      };
    \end{axis}
  \end{tikzpicture}

  ~

  \begin{tabular}{llll} \hline
    Manifold & Volume & Type & As linear combination \\ \hline
    $v3318$ & $6.4506\dotsc$ & geometric & $1\cdot\vol(v3318)+0\cdot\vol(L14n54610)$ \\
    $L14n534610$ & $15.6881\dotsc$ & geometric & $0\cdot\vol(v3318)+1\cdot\vol(L14n54610)$ \\
    $t09825$ & $2.2704\dotsc$ & exotic & $4\cdot\vol(v3318)-1.5\cdot\vol(L14n54610)$ \\
    $v3318$ & $0.8770\dotsc$ & exotic & $5\cdot\vol(v3318)-2\cdot\vol(L14n54610)$ \\
    $10^{3}_{5}$ & $5.5736\dotsc$ & exotic & $-4\cdot\vol(v3318)+2\cdot\vol(L14n54610)$ \\
    $L13n5993$ & $1.9097\dotsc$ & exotic & $-7\cdot\vol(v3318)+3\cdot\vol(L14n54610)$ \\
    $K12n809$ & $3.6638\dotsc$ & exotic & $3\cdot\vol(v3318)-1\cdot\vol(L14n54610)$ \\
    $L13n5993$ & $9.2374\dotsc$ & exotic & $-1\cdot\vol(v3318)+1\cdot\vol(L14n54610)$ \\
    $v2489$ & $1.8319\dotsc$ & exotic & $1.5\cdot\vol(v3318)-0.5\cdot\vol(L14n54610)$ \\
    $t09825$ & $5.0572\dotsc$ & exotic & $2\cdot\vol(v3318)-0.5\cdot\vol(L14n54610)$ \\
    $v3548$ & $2.7868\dotsc$ & exotic & $-2\cdot\vol(v3318)+1\cdot\vol(L14n54610)$ \\
  \end{tabular}
  \caption{$p(x) = x^4+x^2-2x+1$, root $-0.624\dotsc+1.300\dotsc i$,
  prospective basis: $\vol(v3318)$, $\vol(L14n54610)$ with fit ratio 2.}
  \label{fig:v3318-lattice}
\end{figure}

Figure~\ref{fig:v3318-lattice} represents a typical ``almost good''
lattice from our data.  Our fit ratio is 2, which points out that we
are unable to fit the exotic volumes of $v2489$ and $t09825$ (either
of them) completely into our lattice. Nonetheless, this lattice allows
us to conjecture the existence of smaller manifolds which would
``fix'' this lattice. For example, we might suppose that the nearby
exotic volumes of $v2489$ and $t09825$ differ by a generator (or,
perhaps, a multiple of a generator): this would imply the existence of
a manifold with volume $\frac{\vol(v3318)}{2} \approx 3.2253$ (this is
by no means the only volume which would ``fix'' the lattice).

\clearpage

\printbibliography

\end{document}